\documentclass[twoside]{article}

\oddsidemargin=\evensidemargin
\usepackage{amsfonts}

\catcode`\@=11
\long\def\@makefntext#1{
\protect\noindent \hbox to 3.2pt {\hskip-.9pt  
$^{{\eightrm\@thefnmark}}$\hfil}#1\hfill}		

\def\ps@myheadings{\let\@mkboth\@gobbletwo		
\def\@oddhead{\hbox{}
\rightmark\hfil\eightrm\thepage}   
\def\@oddfoot{}\def\@evenhead{\eightrm\thepage\hfil
\leftmark\hbox{}}\def\@evenfoot{}
\def\sectionmark##1{}\def\subsectionmark##1{}}

\catcode`@=11		      		     
\def\ps@plain{\let\@mkboth\@gobbletwo
     \def\@oddhead{}\def\@oddfoot{\eightrm\hfil\thepage
     \hfil}\def\@evenhead{}\let\@evenfoot\@oddfoot}

\newcounter{sectionc}\newcounter{subsectionc}\newcounter{subsubsectionc}
\renewcommand{\section}[1] {\vspace{12pt}\addtocounter{sectionc}{1} 
\setcounter{subsectionc}{0}\setcounter{subsubsectionc}{0}\noindent 
	{\tenbf\thesectionc. #1}\par\vspace{5pt}}
\renewcommand{\subsection}[1] {\vspace{12pt}\addtocounter{subsectionc}{1} 
	\setcounter{subsubsectionc}{0}\noindent 
	{\bf\thesectionc.\thesubsectionc. 
	{\kern1pt \bfit #1}}\par\vspace{5pt}}
\renewcommand{\subsubsection}[1] {\vspace{12pt}
	\addtocounter{subsubsectionc}{1}
	\noindent
	{\tenrm\thesectionc.\thesubsectionc.\thesubsubsectionc.	{\kern1pt 
	\it #1}}\par\vspace{5pt}}

\newcounter{appendixc}
\newcounter{subappendixc}[appendixc]
\newcounter{subsubappendixc}[subappendixc]

\renewcommand{\appendix}[1] {\vspace{12pt}	
	\refstepcounter{appendixc}		
	\setcounter{figure}{0}
	\setcounter{table}{0}
	\setcounter{lemma}{0}
	\setcounter{theorem}{0}
	\setcounter{corollary}{0}
	\setcounter{definition}{0}
	\setcounter{equation}{0}
	\renewcommand{\thefigure}{\Alph{appendixc}.\arabic{figure}}
	\renewcommand{\thetable}{\Alph{appendixc}.\arabic{table}}
	\renewcommand{\theappendixc}{\Alph{appendixc}}
	\renewcommand{\thelemma}{\Alph{appendixc}.\arabic{lemma}}
	\renewcommand{\thetheorem}{\Alph{appendixc}.\arabic{theorem}}
	\renewcommand{\thedefinition}{\Alph{appendixc}.\arabic{definition}}
	\renewcommand{\thecorollary}{\Alph{appendixc}.\arabic{corollary}}
	\renewcommand{\theequation}{\Alph{appendixc}.\arabic{equation}}
	\noindent{\tenbf Appendix \theappendixc #1}\par\vspace{5pt}}

\newcommand{\textlineskip}{\baselineskip=13pt}
\newcommand{\smalllineskip}{\baselineskip=10pt}

\newcommand{\copyrightheading}[1]
	{\vspace*{-2.5cm}\smalllineskip{\flushleft
	{\footnotesize \hspace{.1in}}\\
   	{\footnotesize \hspace{.1in}}\\
         }}

\def\abstracts#1#2#3#4{{
	\centering{\begin{minipage}{4.5in}\footnotesize\baselineskip=10pt
	\centerline{ABSTRACT} 
	\parindent=15pt #1\par 
	\parindent=15pt #2\par
	\parindent=15pt #3\par
	\parindent=15pt #4\par
	\end{minipage}}\par}} 

\def\keywords#1{{ 
	\centering{\begin{minipage}{4.5in}\footnotesize\baselineskip=10pt
	{\footnotesize\it Keywords}\/: #1
	\end{minipage}}\par}}

\renewenvironment{thebibliography}[1]
	{\frenchspacing
	 \ninerm\baselineskip=11pt
	 \begin{list}{[\arabic{enumi}]}
	{\usecounter{enumi}\setlength{\parsep}{0pt}
	 \setlength{\leftmargin 13.7pt}{\rightmargin 0pt} 
	 \setlength{\itemsep}{0pt} \settowidth
	{\labelwidth}{[#1]}\sloppy}}{\end{list}}

\newcounter{itemlistc}
\newcounter{romanlistc}
\newcounter{alphlistc}
\newcounter{arabiclistc}

\newcommand{\fcaption}[1]{
        \refstepcounter{figure}
        \setbox\@tempboxa = \hbox{\footnotesize Fig.~\thefigure. #1}
        \ifdim \wd\@tempboxa > 5in
           {\begin{center}
        \parbox{5in}{\footnotesize\smalllineskip Fig.~\thefigure. #1}
            \end{center}}
        \else
             {\begin{center}
             {\footnotesize Fig.~\thefigure. #1}
              \end{center}}
        \fi}

\newcommand{\tcaption}[1]{
        \refstepcounter{table}
        \setbox\@tempboxa = \hbox{\footnotesize Table~\thetable. #1}
        \ifdim \wd\@tempboxa > 5in
           {\begin{center}
        \parbox{5in}{\footnotesize\smalllineskip Table~\thetable. #1}
            \end{center}}
        \else
             {\begin{center}
             {\footnotesize Table~\thetable. #1}
              \end{center}}
        \fi}

\def\pmb#1{\setbox0=\hbox{#1}
	\kern-.025em\copy0\kern-\wd0
	\kern.05em\copy0\kern-\wd0
	\kern-.025em\raise.0433em\box0}

\def\fnt#1#2{\footnotetext{\kern-.3em
	{$^{\mbox{\scriptsize #1}}$}{#2}}}

\def\runninghead#1#2{\pagestyle{myheadings}
\markboth{{\protect\footnotesize\it{\quad #1}}\hfill}
{\hfill{\protect\footnotesize\it{#2\quad}}}}

\font\tenrm=cmr10
 
\font\tenbf=cmbx10
\font\bfit=cmbxti10 at 10pt
\font\ninerm=cmr9

\font\eightrm=cmr8

\newtheorem{thm}{Theorem}   

\newtheorem{lem}{Lemma}

\newtheorem{cor}{Corollary}
\newtheorem{conj}{Conjecture}
\def\@begintheorem#1#2{\trivlist	
	\item[\hskip\labelsep{\bf #1\ #2.}]} 
\def\@opargbegintheorem#1#2#3{\trivlist
	\item[\hskip\labelsep{\bf #1\ #2\ (#3).}]}

\newenvironment{proof}{\begin{trivlist}
	\item[\noindent]{\it Proof.}}{\quad $\square$\end{trivlist}} 

    	{\setcounter{itemlistc}{0}		
	 \begin{list}{$\bullet$}		
	{\usecounter{itemlistc}			
	 \leftmargin10pt	       
	 \setlength{\parsep}{0pt}
	 \setlength{\itemsep}{0pt}     
	}}{\end{list}}

	{\setcounter{romanlistc}{0}		
	 \begin{list}{$($\roman{romanlistc}$)$}	
	{\usecounter{romanlistc}		
	 \leftmargin18pt 
	 \setlength{\parsep}{0pt}
	 \setlength{\itemsep}{0pt}	
	 \settowidth{\labelwidth}{#1}                          
	}}{\end{list}}

	{\setcounter{enumii}{0}			
	 \begin{list}{$($\alph{enumii}$)$}	
	{\usecounter{enumii}			
	 \leftmargin18pt		
	 \setlength{\parsep}{0pt}
	 \setlength{\itemsep}{0pt}	
	 \settowidth{\labelwidth}{#1}                          
	}}{\end{list}}

\def\qed{\hbox{${\vcenter{\vbox{			
   \hrule height 0.4pt\hbox{\vrule width 0.4pt height 6pt
   \kern5pt\vrule width 0.4pt}\hrule height 0.4pt}}}$}}

\def\theequation{\thesectionc.\arabic{equation}}  

\newcommand{\sph}{\mbox{$\mathbb S^3$}}
\newcommand{\psl}{\mbox{PSL$_2 \mathbb C$}}
\newcommand{\slt}{\mbox{SL$_2 \mathbb C$}}

\usepackage{graphicx}
\begin{document}
\setlength{\textheight}{7.7truein}  

\runninghead{J.\ Hoste \& P.\ D.\ Shanahan}{Commensurability Classes of Twist Knots}

\normalsize\textlineskip
\thispagestyle{empty}
\setcounter{page}{1}

\copyrightheading{}		    

\vspace*{0.88truein}

\centerline{\bf COMMENSURABILITY CLASSES OF TWIST KNOTS}
\vspace*{0.37truein}
\centerline{\footnotesize JIM HOSTE  }
\baselineskip=12pt
\centerline{\footnotesize\it Department of Mathematics}
\baselineskip=10pt
\centerline{\footnotesize\it Pitzer College}
\baselineskip=10pt
\centerline{\footnotesize\it Claremont, CA 91711, USA}
\baselineskip=10pt
\centerline{\footnotesize\it jhoste@pitzer.edu}

\vspace*{10pt}
\centerline{\footnotesize PATRICK D. SHANAHAN}
\baselineskip=12pt
\centerline{\footnotesize\it Department of Mathematics}
\baselineskip=10pt
\centerline{\footnotesize\it Loyola Marymount University}
\baselineskip=10pt
\centerline{\footnotesize\it Los Angeles, CA 90045, USA}
\baselineskip=10pt
\centerline{\footnotesize\it pshanahan@lmu.edu}

\vspace{1in}

\vspace*{0.21truein} 
\abstracts{In this paper we prove that if $M_K$ is the complement of a non-fibered twist knot $K$ in \sph ,
then $M_K$ is not commensurable to a fibered knot complement in a $\mathbb Z/ 2 \mathbb Z$-homology sphere.   To prove
this result we derive a recursive description of the character variety of twist knots and 
then prove that a commensurability criterion developed by D.\ Calegari and N.\ Dunfield is satisfied
for these varieties. In addition, we partially extend our results to a second infinite family of 2-bridge knots.}{}{}{}

\vspace*{10pt}
\keywords{twist knot, commensurable, fibered knot, character variety}


\vspace*{1pt}\textlineskip	
\section{Introduction.}	        
\vspace*{-0.5pt}

Let $M$ be a compact, orientable 3-manifold. Recall that $M$ is {\em fibered} if $M$ is 
homeomorphic to a surface bundle over $\mathbb S^1$, and $M$ is {\em 
virtually fibered} if it has a finite cover which is fibered.  In 1982, Thurston conjectured 
that every finite volume, hyperbolic 3-manifold is virtually fibered \cite{th}. Because there has been little
progress made in resolving the general case of this conjecture, it is natural to focus on
specific classes of manifolds for which the problem is more accessible.   For example, hyperbolic knot complements  in
\sph \ provide an interesting class of 3-manifolds in which to study  Thurston's conjecture. Since every finite cover
of a fibered knot complement is fibered, a sufficient condition for a knot complement to be virtually fibered is for
it to share a finite cover with a fibered knot  complement. Two manifolds which share a common finite cover are called
{\em commensurable}.  

In \cite{cd}, Calegari and Dunfield  investigate the question of when a non-fibered, hyperbolic 
knot complement can be commensurable to a fibered knot complement. Before stating their result we  review some necessary terminology.  If $M$ is a hyperbolic manifold, then let $X(M)$ be the \psl-character variety of $\pi_1(M)$ and let $X_0(M)$ denote the
irreducible component of $X(M)$ containing the character of the discrete faithful representation. 
Furthermore, call
$M$ {\em generic} if it is not arithmetic and its commensurator orbifold has a flexible cusp. The latter condition holds if the cusp coefficient of $M$ is not in either of the fields $\mathbb Q(i)$ or $\mathbb Q (\sqrt{-3})$. For example, all
non-fibered twist knots in
\sph \ have generic hyperbolic complements. To see this, first note that Reid proved in \cite{re} that the only arithmetic knot in \sph \
is the figure-eight knot (which is fibered). Moreover, the computation of the degree of the trace field of twist knots
in
\cite{hsi}, together with the fact that the trace field and the cusp field of these knots are identical (see
\cite{nr}), implies that the cusp coefficient of a non-fibered twist knot has degree greater than 2 over $\mathbb Q$. We
are now prepared to state Calegari and Dunfield's theorem.

\begin{thm}[Calegari and Dunfield] Let $M$ be a generic hyperbolic knot complement in a $\mathbb Z/2 \mathbb Z$-homology
sphere. Suppose that $X_0(M)$ contains the character of a non-integral reducible representation. Then $M$ is not
commensurable to a fibered knot complement in a $\mathbb Z/2 \mathbb Z$-homology sphere.
\label{crit}
\end{thm}

\begin{figure}[htbp] 
\vspace*{13pt}
\centerline{\includegraphics[scale=.75]{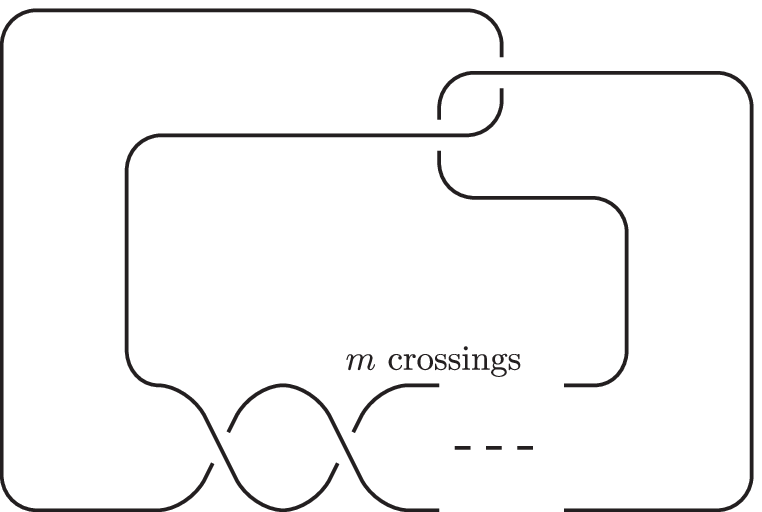}}
\vspace*{13pt}
\fcaption{The twist knot $K_m$.}
\label{figone}
\end{figure}

Using this theorem and results of Hilden, Lozano, and Montesinos \cite{hlm}, Calegari and Dunfield show that all non-fibered 2-bridge knots
$K_{p/q}$ with $0<p<q<40$ have complements which are not commensurable with a fibered knot complement. 
They remark that it would have been nicer to establish this result for all non-fibered 2-bridge
knots. In this article we extend Calegari and Dunfield's result by proving that the infinite family of non-fibered
{\em twist} knots $K_m$ satisfy the hypotheses of Theorem~\ref{crit} where $K_m$ is pictured in Figure~\ref{figone}.
This gives the following theorem.

\begin{thm} If $M_{K_m}$ is the complement of the non-fibered twist knot $K_m$ in $\mathbb S^3$, then $M_{K_m}$ is not
commensurable to a fibered knot complement in a $\mathbb Z/ 2 \mathbb Z$-homology sphere.
\label{main}
\end{thm}

Theorem~\ref{main} does not imply that non-fibered twist knot complements are
not virtually fibered. In fact, Leininger \cite{le}  proved that twist knot complements 
are virtually fibered,  and subsequently, Walsh \cite{wa} extended this result to 2-bridge knots.

The paper will proceed as follows. In Section~2, we derive a recursive description for the \psl-character variety of twist knot complements. Our procedure extends to other infinite families of 2-bridge knots described in Section~4. In Section~3 we turn our attention to proving Theorem~\ref{main} by verifying that the character varieties described in
Section~2 satisfy the hypotheses of Theorem~\ref{crit}. The key ingredient of the proof is establishing that the character varieties of twist knots are $\mathbb Q$-irreducible. Finally, in Section~4, we partially extend these results to a second class of 2-bridge knots.

\section{The Character Variety of Twist Knots}

In this section we develop a recursive description of the character variety of twist knots. The diagram of $K_m$ 
shown in Figure~\ref{figone} has $|m|+2$ crossings: 2 in the ``clasp'' and $|m|$ in the ``twist.'' Our convention is
that the twists are right-handed if
$m>0$ and left-handed if $m<0$.  Consider the complement of the knot $K_m$ in $\mathbb S^3$. After
reflecting in the plane of the diagram and applying a twist to change the clasp back to a right-handed one, we see
that the complement of the twist knot $K_m$ is homeomorphic to the complement of $K_{1-m}$.
Therefore, from this point on we shall assume that $m=2n$ and we shall let $M_n$ denote the
complement of $K_{2n}$ in \sph.  It is well known that  the only non-hyperbolic twist knots are the unknot and the
trefoil. Thus,
$M_n$ is  a hyperbolic manifold for $n
\neq 0, 1$. Recall also that the only fibered twist knots are the unknot, trefoil, and the figure-eight, hence
$M_n$ is non-fibered if $n \neq 0,1,-1$. 

The fundamental group of $M_{n}$ has a presentation of the form
\begin{equation}
\pi_1(M_n)  =  \left< a, b\ | \ a w^n = w^n b \right>
\label{pres}
\end{equation}
where $a$ and $b$ are meridians, and $w=(ba^{-1}b^{-1}a)^{-1}$. A slightly different but isomorphic presentation is
derived in Proposition~1 of \cite{hsii}.  We have chosen to use the present form in order to agree with the
presentation used in \cite{hlm}.

A {\em representation} of $\pi_1(M_n)$ in \slt \ is a group homomorphism $\rho : \pi_1(M_n)
\rightarrow \slt$. A representation is called {\em reducible} if there exists a one-dimensional eigenspace of the image, otherwise it is called {\em irreducible}. Equivalently, $\rho$ is reducible if all matrices in the image $\rho(\pi_1(M_n))$ 
can be simultaneously conjugated to be upper triangular. A representation is {\em abelian} if the image is
an abelian subgroup of \slt. Every abelian
representation is reducible. On the other hand, there do exist reducible, non-abelian representations.

If we let $R(M_n)$ denote the set of all representations of $\pi_1(M_n)$ in \slt, then $R(M_n)$ has the
structure of an affine algebraic set in $\mathbb C^{8}$ determined by six polynomial equations: four coming from the single relation in (\ref{pres}) and two coming from the requirement that
$\mbox{det}(\rho(a))=\mbox{det}(\rho(b))=1$. In order to apply Theorem~\ref{crit}, we need to determine the \psl-character variety of $\pi_1(M_n)$. The {\em character} of a representation $\rho \in R(M_n)$ is
the function $\chi_\rho:\pi_1(M_n) \rightarrow \mathbb C$ defined by $\chi_\rho(g)=\mbox{trace}(\rho(g))$. The
set of all characters forms an affine algebraic set $X(M_n)$ which is called the {\em character variety}.
Two representations $\rho, \rho' \in R(M_n)$ are called {\em conjugate} if there is a matrix $A$ in \slt \ such that
$\rho(g) = A \rho'(g) A^{-1}$ for all $g \in \pi_1(M_n)$. Conjugacy defines an equivalence relation on
$R(M_n)$, and since trace is invariant under conjugation, the character variety may be computed from the
quotient of $R(M_n)$ under this relation. Moreover, because Theorem~\ref{crit} is concerned only with
the  component $X_0(M_n)$ which contains the character of the discrete, faithful representation and because this
representation is irreducible (hence, non-abelian) we will initially restrict our attention to the set $\hat R(M_n)$ of
conjugacy classes of non-abelian representations.

Now assume that $\rho \in R(M_n)$ is a non-abelian representation. Then $\rho$ may be conjugated so that:
\begin{equation}
\rho(a) = \left( \begin{array}{cc} m & 1 \\ 0 & 1/m \end{array} \right) \ \ \ \mbox{and}\ \ \ 
\rho(b) = \left( \begin{array}{cc} m & 0 \\ -q & 1/m \end{array} \right).
\label{images}
\end{equation}
Here the choice of $-q$ as opposed to $q$ is entirely arbitrary.  To avoid cumbersome notation we identify
$a$ with $\rho(a)$, $w^n$ with $\rho(w^n)$, and so on. It is easy to compute the eigenvalues and eigenvectors of $a$ and $b$. One can then show that $a$ and $b$ share a common one-dimensional eigenspace if and only if $q=0$ or $q=(m-1/m)^2$.

Notice that an ordered pair $(m,q) \in \mathbb C^2$
corresponds to a non-abelian representation if and only if the relation $a w^n = w^n b$ from (\ref{pres}) holds in
\slt . If we view $m$ and $q$ as variables, then this relation leads to four polynomial equations in $m$, $1/m$, and
$q$ which must be satisfied. In \cite{ri}, Riley proves that the four
equations reduce to the single equation:
\begin{equation}
R_n = (m-1/m) w^n_{12} + w^n_{22}=0,
\label{rieq}
\end{equation}
where $w^n_{ij}$ denotes the $(i,j)$-entry of the matrix $w^n=\rho(w^n)$. Thus, the algebraic set $\hat R(M_n)$ is
defined by the polynomial $R_n$ given in (\ref{rieq}). 

In \cite{hsii} we derive a recursive formula for $R_n$ as follows. From the Cayley-Hamilton Theorem and the fact that
$\mbox{det}(w)=1$ we have that
\begin{equation}
w^{n+1} - T \, w^{n} + w^{n-1}=0,
\label{wrecur}
\end{equation}
where $T = \mbox{trace}(w)$. Thus, the entries of $w^n$ satisfy the same recursion, and therefore, from
(\ref{rieq}) we obtain:
\begin{equation}
R_{n+1} - T \, R_{n} + R_{n-1}=0,
\label{Rrec}
\end{equation}
where
$$T = \mbox{trace}((ba^{-1}b^{-1}a)^{-1})= 2+(2-m^2-1/m^2)q+q^2$$
and with initial conditions
$$\begin{array}{lcl}
R_0    & = & 1, \\
R_1    & = & -1 +m^2 + 1/m^2 - q.
\end{array}$$

Observe that if $q=(m-1/m)^2$, then $R_1=1$ and $T=2$ implying that $R_n=1$ for all $n$.  Thus, if $\rho$ is a representation then $q \neq (m-1/m)^2$.  In other words, $\rho$ is a non-abelian reducible representation if and only if $q=0$ and $m \neq \pm 1$.

In \cite{hlm}, Hilden, Lozano, and Montesinos show that the \psl-character variety is parameterized by the variables $x=\mbox{trace}(\rho(a^2))$ and $y=\mbox{trace}(\rho(ab))$.  From (\ref{images}) we see that $x$ and $y $ have the following values
\begin{eqnarray} \label{xzvalues}
x & = & m^2 + 1/m^2, \\ \nonumber
z & = & m^2 + 1/m^2-q. 
\end{eqnarray}
By changing variables from $m$ and $q$ to $x$ and $z$ in our polynomial $R_n$ we obtain a polynomial that defines the \psl-character variety $X(M_n)$.

\begin{lem}
Let $r_n(x,z) \in \mathbb Z[x,z]$ be defined recursively by
\begin{equation}
r_{n+1}(x,z) - t(x,z) \, r_{n}(x,z) + r_{n-1}(x,z)=0,
\label{rrec}
\end{equation}
where
$$t(x,z) = 2+2x-2z-xz+z^2$$
and with initial conditions
$$\begin{array}{lcl}
r_0(x,z)    & = & 1, \\
r_1(x,z)    & = & z-1.
\end{array}$$
Then $(x-z)r_n(x,z)$ is the defining polynomial of the \psl-character variety $X(M_n)$. Moreover, the defining
polynomial of $X_0(M_n)$ is a factor of $r_n(x,z)$. 
\label{rlem}
\end{lem}

\begin{proof} If $\rho$ is an abelian representation then $\rho(a)=\rho(b)$ because $a$ and $b$ are meridians, and so $x=z$. Conversely, for every matrix $A \in \slt$, the mapping $\rho$ defined by $\rho(a)=\rho(b)=A$ is an abelian representation. Therefore, the polynomial $x-z$ defines the characters of all abelian representations.  For non-abelian
representations, it suffices to verify that the substitution of $x=m^2+1/m^2$ and $z=m^2+1/m^2-q$ gives the 
recursive description of $R_n$. The final statement of the lemma follows from the fact that the discrete faithful
representation is non-abelian, and so its character does not lie on the curve defined by $x-z$.
\end{proof}

We remark that the polynomials $r_n$ defined
above agree with the polynomials
$r[(-4n+1)/(-2n+1)]$ if
$n<0$ and
$r[(4n-1)/(2n-1)]$ if $n>0$ defined in \cite{hlm}. Our recursive description, however, is significantly different.

\section{Calegari and Dunfield's Criterion}

In this section we show that the character varieties of non-fibered twist knots satisfy the hypotheses of Theorem~\ref{crit}. As mentioned in the introduction, all non-fibered twist knot complements are generic. Thus, it remains to show that $X_0(M_n)$ contains the character of a non-integral reducible representation. Recall that a representation $\rho : \pi_1(M_n) \rightarrow \slt$ is called {\em integral} if trace$(\rho(g))$ is an algebraic integer for all $g \in \pi_1(M_n)$.  Otherwise, $\rho$ is called {\em non-integral}.

By Lemma~\ref{rlem} if $X_0(M_n)$ contains the character of a non-integral reducible representation then this character must lie on the curve defined by $r_n(x,z)$. Furthermore, as stated earlier, a non-abelian representation $\rho$ as defined by  (\ref{images}) is reducible if and only if $q=0$ and $m \neq \pm 1$.  For such a representation it follows from (\ref{xzvalues}) that $x=z$. Thus, if $X_0(M_n)$ contains the character of a non-integral reducible representation, then this representation corresponds to a root of the polynomial $r_n(x,x)$. These polynomials are easily determined by setting $z=x$ in the formula for $r_n$ from Lemma~\ref{rlem}. This gives the following corollary.

\begin{cor} $r_n(x,x)= n\, x -(2n-1).$
\end{cor}

Clearly, if $n \neq 0, \pm 1$, then $x= (2n-1)/n$ is a root of $r_n(x,x)$ which is not an algebraic integer.
However, this alone does not imply that the the twist knots satisfy the hypotheses of Theorem~\ref{crit}. While the value
$x= (2n-1)/n$ for $n \neq 0, \pm 1$ does correspond to the character of a non-integral  reducible representation in
$X(M_n)$, it is not yet clear that this character lies in the geometric component $X_0(M_n)$.  

The requirement that the character of a non-integral reducible representation lie in the geometric component is
essential in Calegari and Dunfield's proof of Theorem~\ref{crit}. In particular, they show that if
$M_1$ and $M_2$ are commensurable, generic, 1-cusped manifolds, then there is a natural birational isomorphism between
$X_0(M_1)$ and $X_0(M_2)$ (see, for example \cite{lr}). They use this isomorphism to show that if
$X_0(M_1)$ contains the character of a non-integral reducible representation, then so does $X_0(M_2)$. However, if
$M_2$ is fibered, then it is known that $X_0(M_2)$ cannot contain such a character. This gives the criterion of
Theorem~\ref{crit}. Additionally, in Remark~7.2 of \cite{cd}, Calegari and Dunfield point out that the maps used in their proof are all
defined over $\mathbb Q$ as opposed to $\mathbb C$. Therefore, the hypotheses of Theorem~\ref{crit} can be
weakened to having a non-integral reducible representation in the $\mathbb Q$-irreducible component of $X(M_1)$
containing $X_0(M_1)$. Returning to the case of twist knots, if we show that $r_n(x,z)$ is $\mathbb
Q$-irreducible, then the curve defined by $r_n(x,z)$ is the $\mathbb Q$-irreducible component of $X(M_n)$ containing $X_0(M_n)$. Therefore, the non-integral reducible representation corresponding to $x=(2n-1)/n$ would satisfy the
(weakened) hypotheses of Theorem~\ref{crit}. 

\begin{lem} For all $n \in \mathbb Z$, the polynomial $r_n(x,z)$ is $\mathbb Z$-irreducible.
\label{rirred}
\end{lem}

\begin{proof} This is clear for the two non-hyperbolic values $n = 0,1$. So, for $n \neq 0,1$ assume by way of
contradiction that $r_n(x,z)$ factors as $r_n(x,z)=f_n(x,z)g_n(x,z)$ with neither $f_n$ nor $g_n$ a
constant. Consider for the moment the case $n<0$. From the recursion in Lemma~\ref{rlem}, we see that the total degree
of $r_n(x,z)$ is $-2n$. Therefore, every term in both $f_n$ and $g_n$ has total degree strictly less than $-2n$. Now
if
$\rho_0$ is the discrete faithful representation of $\pi_1(M_n)$, then we must have:
$$
\rho_0(a) = \left( \begin{array}{cc} 1 & 1 \\ 0 & 1 \end{array} \right) \ \ \ \mbox{and}\ \ \ 
\rho_0(b) = \left( \begin{array}{cc} 1 & 0 \\ -q_0 & 1 \end{array} \right)
$$
where $q_0$ is a particular root of the polynomial $R_n(1,q)$. In Section~3 of \cite{hsi}, we show that the polynomial
$R_n(1,q)$ (which is denoted by $\Phi_{-n}(q)$ in that paper) is $\mathbb Z$-irreducible and of degree $-2n$. Thus,
$q_0$ has degree $-2n$ over $\mathbb Q$. The representation $\rho_0$ given by the zero $(1,q_0)$ of $R_n(m,q)$ has a character given by the zero $(2,2-q_0)$ of $r_n(x,z)$. That is, $q_0$ is a root of the polynomial $r_n(2,2-q)=f_n(2,2-q)g_n(2,2-q)$ and so either
$f_n(2,2-q_0)=0$ or $g_n(2,2-q_0)=0$.  However, each of $f_n(2,2-q)$ and $g_n(2,2-q)$ have degree
strictly less than $-2n$ in $q$. Thus, the degree of $q_0$ over $\mathbb Q$  is less than $-2n$, a contradiction. The
proof for
$n>0$ is similar.
\end{proof}

This completes the proof of Theorem~\ref{main}.  The strategy of the proof can be summarized in the following corollary of Calegari and Dunfield's Theorem.
\begin{thm} Let $M$ be a generic, hyperbolic, non-fibered 2-bridge knot complement in $\mathbb S^3$.  If the total degree of $r_M(x,z)$ is equal to the degree of $r_M(2,2-q)$, if $r_M(x,x)$ is not monic, and if 
$r_M(2,2-q)$ is $\mathbb Z$-irreducible, then $M$ is not commensurable to a fibered knot complement in a $\mathbb Z/ 2 \mathbb Z$-homology sphere.
\label{newcrit}
\end{thm}

\section{Further Results}

\begin{figure}[htbp] 
\vspace*{13pt}
\centerline{\includegraphics[scale=.75]{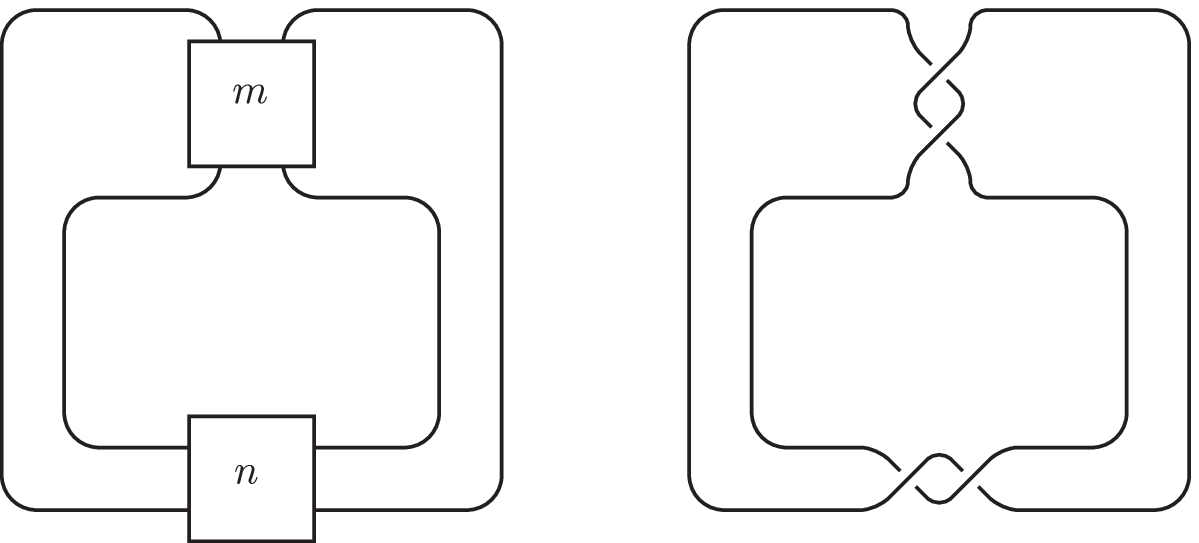}}
\vspace*{13pt}
\fcaption{$J(m,n)$ and the figure eight knot $J(2,-2)$.}
\label{jmn}
\end{figure}

The twist knots are part of a slightly more general family of knots 
$J(m,n)$ defined by Figure~\ref{jmn}, where $m$ and $n$ are the
number of right-handed crossings contained in each box. Note that $J(m,n)$
is a knot precisely when $mn$ is even, and a
two-component link otherwise.  For example, $J(-1,2)$ is the
right-handed trefoil, $J(2,-2)$ is the figure eight knot, and $J(-1,1)$ is
the Hopf link. If $m=1$ we obtain a torus knot, while $m=2$ gives a twist
knot.

Clearly, $J(m,n)$ is symmetric in $m$ and $n$. If $J(m,n)$ is a knot, then its
orientation is immaterial, as there exists an obvious rotation of $\mathbb S^3$
carrying the knot to its reverse.  Furthermore, $J(-m,-n)$ is the mirror image
of $J(m,n)$. Thus, we may consider, for example, only those $m$ and $n$ for which 
$m>0$ and $n$ is even. 

By Proposition~1 of \cite{hsii} the fundamental group of $J(m,2n)$ has a presentation of the form
$$\pi_1(\mathbb S^3 - J(m,2n))=<a, b\ |\ aw_m^n = w_m^n b>$$
where $w_m$ is a word in $a$ and $b$ given by a formula depending only on $m$.  This presentation  allows us to repeat the calculations of Section~2 for any family $\{ J(m,2n) \}_{n=-\infty}^{\infty}$ with $m$ fixed.

\begin{figure}[htbp] 
\vspace*{13pt}
\centerline{\includegraphics[scale=.60]{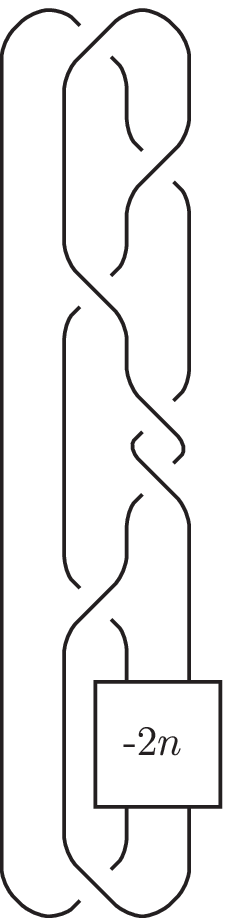}}
\vspace*{13pt}
\fcaption{$J(-3,-2n)$.}
\label{jplatt}
\end{figure}

For the remainder of this paper we will focus on extending the results of Section~3 to the knots $J(3,2n)$.  We will first show that these knots are fibered if and only if $n \ge 0$.  Its not hard to see that the mirror image of $J(3,2n)$ can be redrawn as in Figure~\ref{jplatt}. As illustrated in \cite{ht} this diagram determines the continued fraction
$$[1,-2,-2n]=\frac{1}{1-\frac{1}{-2-\frac{1}{-2n}}}=\frac{4n-1}{6n-1}.$$
Recall that two 2-bridge knots with fractions $p/q$ and $p'/q'$ are ambient isotopic if and only if $q'=q$ and $p' \equiv p^{\pm1}  (\mbox{mod } q)$. Therefore, the fraction $(4n-1)/(6n-1)$ gives the same 2-bridge knot as the fraction $(6n-4)/(6n-1)$.  If $n > 0$, then
$$\frac{6n-4}{6n-1} = [2, 2, \dots , 2, -2n]$$
where the continued fraction has $2n$ entries. From this description we see that the knot is a band-connected sum of Hopf links and therefore is fibered \cite{ha}.  

On the other hand, for $n<0$ a routine calculation shows that the Alexander polynomial of $J(3,2n)$ is
$$\Delta(t)=2-3t+3t^2- \dots -3t^{-2n-1}+2t^{-2n}.$$
Since these Alexander polynomials are not monic, these knots are not fibered.
We would like to show that none of the knots $J(3,2n)$ for $n<0$ have a complement which  is commensurable to a fibered knot complement in a $\mathbb Z/2 \mathbb Z$-homology sphere.
First note that all of these knots are hyperbolic;  by Theorem~1 of \cite{ht} none are satellites and furthermore none are torus knots.  

In order to apply Theorem~\ref{newcrit} we need to show that the complements of $J(3,2n)$ for $n<0$ are generic.  Our strategy is to first show that the cusp and trace fields coincide for these knots. This is conjectured by Callahan and Reid  to be true for all hyperbolic knots except the pair of dodecahedral knots \cite{cr}.
As already mentioned, Nuemann and Reid proved that the cusp and trace fields coincide for twist knots, and Callahan and Reid remark that the same is true for all hyperbolic knots with eight or less crossings.  

\begin{lem} For $n \neq 0$, the cusp field and trace field of $J(3,2n)$ are identical.
\label{cusp=trace}
\end{lem}

\begin{proof}  In \cite{hsii} we show that the fundamental group of $J(3,2n)$ has the presentation
$$\pi_1(\mathbb S^3 - J(3,2n))=<a, b\ |\ aw^n = w^n b>,$$
where $w=ab^{-1}aba^{-1}b$. Suppose $\rho$ is a parabolic representation defined by (\ref{images}) with $m=1$, then 
$$w= \left( \begin{array}{cc} 1-2q-3q^2-q^3 & 1+2q+q^2 \\ -q(1+2q+q^2) & 1+q+q^2 \end{array} \right).$$
It is easy to verify for $n=0,1$ that 
\begin{equation}
w_{12}^n+q w_{21}^n = 0,\ \mbox{and}
\label{w12rel}
\end{equation}
\begin{equation} 
(1+q)w_{11}^n+q(3+q)w_{12}^n-(1+q)w_{22}^n = 0,
\label{wentryrelation}
\end{equation}
where $w_{ij}^n$ is the $(i,j)$-entry of $w^n$. Now using (\ref{wrecur}) it follows that (\ref{w12rel}) and (\ref{wentryrelation}) are true for all $n$.

In \cite{hsii} we show that the preferred longitude of $J(3,2n)$ is represented by
$$\lambda_n = \left( \begin{array}{cc} w_{11}^n w_{22}^n - q (w_{12}^n)^2 & 2 w_{11}^n w_{12}^n \\
-2q w_{11}^n w_{22}^n & w_{11}^n w_{22}^n-q (w_{12}^n)^2 \end{array} \right) \left( \begin{array}{cc} 1 & -2n \\
0 & 1 \end{array} \right).$$
If $q=q_n$ determines the discrete, faithful representation of $\pi_1(\mathbb S^3 - J(3,2n))$, then it follows from (\ref{rieq}) that
$w_{22}^n(q_n)=0$. Furthermore, $q_n (w_{12}^n(q_n))^2=1$ since det$(w^n)=1$.  Therefore, 
$$\lambda_n(q_n) =  \left( \begin{array}{cc} -1 & 2w_{11}^n(q_n)w_{12}^n(q_n)+2n \\
0 & -1 \end{array} \right),$$
and so the cusp field is $\mathbb Q(2w_{11}^n(q_n)w_{12}^n(q_n)+2n)=\mathbb Q(2w_{11}^n(q_n)w_{12}^n(q_n))= \mathbb Q(\alpha_n)$ where $\alpha_n = 2w_{11}^n(q_n)w_{12}^n(q_n)$.  After multiplying (\ref{wentryrelation}) by $2 w_{12}^n(q_n)$ we obtain
$$\alpha_n + \alpha_n q_n +2 q_n +6 =0.$$
From this relationship it follows that the cusp field $\mathbb Q(\alpha_n)$ is the same as the field $\mathbb Q(q_n)$. Since the trace of the word $ab$ is $2-q_n$, it follows that $\mathbb Q(q_n)$ is the trace field of $J(3,2n)$, and therefore the cusp field and trace fields are identical. \end{proof}

\begin{thm} For $-33<n<0$, the complement of $J(3,2n)$ is not commensurable to a fibered knot in a $\mathbb Z/2 \mathbb Z$-homology sphere.
\end{thm}

\begin{proof} Following the procedure in Section~2, the character variety $X(\mathbb S^3 - J(3,2n))$ is determined by $(x-z)r_n(x,z)$ where $r_n(x,z)$ is defined recursively by (\ref{rrec})
with
\begin{eqnarray*}
t(x,z) & = & -4 x - 2 x^2 + 5 z + 6 x z + x^2 z - 4 z^2 - 2 x z^2 + z^3, \\
r_0(x,z) & = & 1, \\
r_1(x,z) & = & 3 + 2 x - 3 z - x z + z^2.
\end{eqnarray*}
Using induction it is easy to verify for $n<0$ that the total degree of both $r_n(x,z)$ and $r_n(2,2-q)$ equals $-3n$, as well as the fact that the leading term of $r_n(x,x)$ is $2x^{-n}$.  
By Theorem~\ref{newcrit} it remains to show that $\mathbb S^3 - J(3,2n)$ is generic and $r_n(2,2-q)$ is $\mathbb Z$-irreducible.  Irreduciblity was verified using {\em Mathematica} for $-33<n<0$. Moreover, since the cusp and trace fields are identical for these knots by Lemma~\ref{cusp=trace}, and the degree of the trace field is $-3n$, it follows that they are generic.
\end{proof}

Our results suggest the following conjecture.

\begin{conj} The complement $\mathbb S^3 - J(m,n)$ is commensurable to a fibered knot in $\mathbb Z/ 2 \mathbb Z$-homology sphere if and only if $J(m,n)$ is fibered.
\end{conj}

\end{document}